\documentclass[journal]{IEEEtran}
\usepackage{amsmath,amsfonts}
\usepackage{algorithmic}
\usepackage{algorithm}
\usepackage{array}
\usepackage{subfig}
\usepackage{textcomp}
\usepackage{stfloats}
\usepackage{url}
\usepackage{verbatim}
\usepackage{graphicx}
\usepackage{epstopdf}
\usepackage{cite}
\usepackage{xcolor} 
\usepackage[normalem]{ulem}
\def\BibTeX{{\rm B\kern-.05em{\sc i\kern-.025em b}\kern-.08em
    T\kern-.1667em\lower.7ex\hbox{E}\kern-.125emX}}
\usepackage{balance}

\begin{document}

\title{\huge Single-Frequency Self-Alignment RF Resonant Beam for Information and Power Transfer}

\author{
    \normalsize Qingwei Jiang, Mingqing Liu, Mengyuan Xu, Wen Fang, Mingliang Xiong, Qingwen Liu,~\IEEEmembership{Senior Member,~IEEE}
    and Shengli Zhou,~\IEEEmembership{Fellow,~IEEE}

    \thanks{
        This work was supported 
        in part by the National Key Research and Development Program of China under No.2022YFA1004700; 
        in part by the National Natural Science Foundation of China under Grant 62305019, Grant 62301308, Grant 62371342, and Grant 62071334; 
        in part by the Shanghai Municipal Science and Technology Major Project under Grant 2021SHZDZX0100; 
        in part by the Shanghai Municipal Commission of Science and Technology Project under Grant 19511132101; 
        in part by the China National Postdoctoral Program for Innovative Talents under Grant BX20230223; 
        in part by the China Postdoctoral Science Foundation under Grant 2023M732262; 
        in part by Aeronautical Science Foundation of China under Grant 20230007038001; 
        and in part by the Fundamental Research Funds for the Central Universities under Grant 22120210543. 
    }

    \thanks{
        Q. Jiang, M. Xu, M. Xiong and Q. Liu
        are with the College of Computer Science and Technology, Tongji University, Shanghai 201804, China
        (e-mail: jiangqw@tongji.edu.cn, xumy@tongji.edu.cn, mlx@tongji.edu.cn, qliu@tongji.edu.cn).
    }

    \thanks{
        M. Liu and W. Fang 
        are with the College of Electronics and Information Engineering, Tongji University, Shanghai 201804, China
        (e-mail: clare@tongji.edu.cn, wen.fang@tongji.edu.cn).
    }
    
    \thanks{
        S. Zhou is with Department of Electrical and Computer Engineering, University of Connecticut, Storrs, CT 06250, USA (e-mail: shengli.zhou@uconn.edu).
    }
}

\maketitle

\begin{abstract}
    Due to power attenuation, improving transmission efficiency in the radio-frequency (RF) band remains a significant challenge, which hinders advancements in various fields of the Internet of Things (IoT), such as wireless power transfer (WPT) and wireless communication. Array design and retro-directive beamforming (RD-BF) techniques offer simple and effective ways to enhance transmission efficiency. However, when the target is an array or in the near field, the RD-BF system (RD-BFS) cannot radiate more energy to the target due to phase irregularities in the target region, resulting in challenges in achieving higher efficiency. To address this issue, we propose the RF-based resonant beam system (RF-RBS), which adaptively optimizes phase and power distribution between transmitting and receiving arrays by leveraging the resonance mechanism to achieve higher transmission efficiency. We analyze the system structure and develop an analytical model to evaluate power flow and resonance establishment. Numerical analysis demonstrates that the proposed RF-RBS achieves self-alignment without beam control and provides higher transmission efficiency compared to RD-BFS, with improvements of up to 16\%. This self-alignment capability allows the system to effectively transfer power and information across varying distances and offsets. The numerical results indicate the capability to transmit watt-level power and achieve 21 bps/Hz of downlink spectral efficiency in indoor settings, highlighting the advantages of RF-RBS in information and power transfer for mobile applications.
\end{abstract}


\section{Introduction}
\label{sec:intro}

\begin{figure}
	\centering
	\includegraphics[width=3.3in]{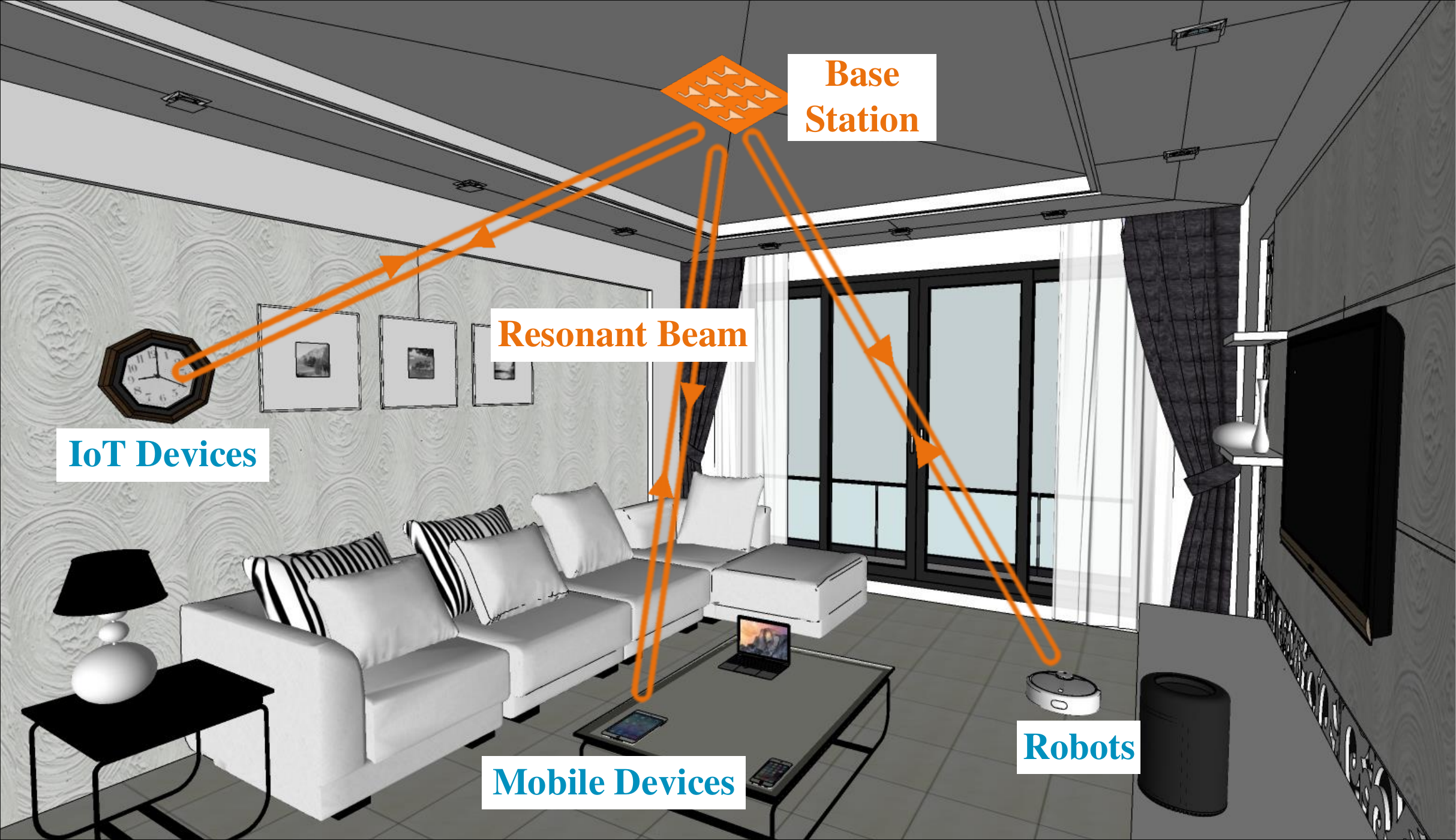}
	\caption{A typical application scenario of RF-RBS in IoT (e.g., WPT, wireless communication).}
	\label{fig: Room scene}
\end{figure}

\IEEEPARstart{T}{he}
advancement of Internet of Things (IoT) technologies and the widespread application of IoT devices \cite{lee2015internet} have facilitated pervasive sensing and communication, playing a significant role in intelligent logistics \cite{zhang2018framework}, smart homes \cite{li2019novel}, and smart cities \cite{kolozali2018observing}. However, the limited battery life of IoT devices challenges their long-term functionality \cite{liu2017green, donatiello2018modeling}. Radio-frequency (RF) wireless power transfer (WPT) technology has emerged as a promising solution, offering a continuous and uninterrupted power supply for IoT devices and the potential to enhance their sensing and communication capabilities \cite{kellogg2014wi}. Nevertheless, RF WPT faces challenges with low end-to-end transmission efficiency due to significant power attenuation of electromagnetic (EM) waves, resulting in high input power demands and a constrained coverage range \cite{visserRFEnergyHarvesting2013, choi2018distributed}. Therefore, improving power transmission efficiency in the RF band remains a critical challenge. 

Array-based designs have shown promise in achieving efficient transmission but still pose a critical concern, commonly referred to as the beam shaping problem \cite{massaArrayDesignsLongDistance2013, StrassnerMicrowavePowerTransmission2013}. The objective is to optimize the phase and power distribution within the transmitting array to maximize efficiency. Retro-directive antenna array (RAA) and beamforming (BF) techniques offer a straightforward and efficient approach for rapid alignment and power transmission. Its receiver only needs to send a pilot signal, enabling the transmitting array to determine the optimal phase distribution and focus EM waves on the receiver \cite{miyamoto2002retrodirective, miyamoto2003digital, nepa2017near, kim2023curved}. To realize this, the Van Atta array \cite{van1959electromagnetic} and the phase-conjugating array \cite{PonRetrodirective1964} are the most common architectures. These techniques, due to their advantages, have been applied in various contexts such as WPT \cite{ettorre2017radiative, re2019circularly}, wireless communications \cite{miyamoto2002retrodirective}, radar \cite{toh2000retrodirective}, RFID systems \cite{chiu2006retrodirective}, passive sensors \cite{islam2017designing}, and interference rejection \cite{goshi2005secure, chun2011interleaved}. Given these applications, it is intuitive that a larger receiver may capture more EM power than a small one. However, when using an antenna array at the receiver to further enhance transmission efficiency, the Retro-Directive Beamforming System (RD-BFS) may not yield higher efficiency, as phase irregularities in the target zone degrade power reception \cite{chou2020conformal}.

Resonance mechanisms exhibit the characteristics of high-efficiency power and information transfer, where two objects exchange energy, as the resonant field is formed by the superposition of in-phase waves, canceling out other phases \cite{aokiObservationStrongCoupling2006}. The mechanisms are applicable across various physical systems, such as acoustic, EM, astronomic, and nuclear. Specifically, in the low or medium-frequency bands, the magnetically coupled resonance WPT (MCR-WPT) system \cite{kurs2007wireless, sample2010analysis} establishes resonant coupling between the transmitter and receiver, resulting in improved power transfer and extended transmission distance. In the lightwave band, lasers \cite{hodgson1997optical} exploit resonance to improve transmission efficiency over natural light sources. Significantly, the optical resonant beam charging (RBC) system \cite{liuChargingSmartphoneAir2022} further establishes a long-range and self-alignment WPT system, transferring over $5$-W optical power across a $2$-m distance with negligible diffraction loss. This optical RBC system also finds applications in communication \cite{xiongRetroReflective2021} and positioning \cite{liuSimultaneousLocalization2023}. However, due to the large sizes of resonance coils in MCR-WPT \cite{rong2021critical} and the low electro-optic and photoelectric conversion efficiency in optical systems, the applications in IoT are limited. In the RF band, masers are microwave generators based on stimulated radiation and resonance \cite{GordonTheMaser1955, OxborrowRoommaser2012}. However, masers' integrated resonant metal cavity limits the application as long-range systems. 

To address these limitations and enhance high-efficiency transmission within the RF band in the near field, we propose the RF-based resonant beam system (RF-RBS). This system integrates RAAs and the resonance mechanism to facilitate self-alignment and achieve high-efficiency transmission for wireless power and information transfer. 
Figure~\ref{fig: Room scene} illustrates a smart home scenario where the RF-RBS is applied, where the base station (BS) is mounted on the ceiling, allowing IoT and mobile devices to be wirelessly charged and communicate simultaneously. The RF-RBS enhances power transfer efficiency and communication reliability while improving radiation safety and reducing EM interference, contributing to the advancement of IoT technologies. 
Based on this conceptual framework, a simultaneous wireless information and power transfer (SWIPT) system is proposed with a focus on the design and performance analysis of a dual-frequency communication scheme \cite{Xia2024Millimeterwave}. However, the mechanism for achieving self-alignment and enhancing efficiency in RF-RBS remains unclear, and an analytical model for evaluating system performance under the three-dimensional (3D) movement of the mobile target (MT) is lacking. 
Additionally, a single-frequency design is preferred over the dual-frequency design for resonance establishment, as it reduces antenna configuration and hardware complexity. Thus, this work provides a detailed system design and develops models to elucidate the underlying mechanisms and analyze the performance of RF-RBS operating in a single frequency, particularly in terms of efficiency and mobility. Simulations incorporating the characteristics of practical components are conducted to validate the feasibility of the proposed RF-RBS approach. 

The contributions of this work are as follows:

\begin{itemize}
    \item We propose a novel single-frequency RF-RBS architecture with low hardware complexity, leveraging the resonance mechanism and RAAs at both the BS and UD to achieve energy-concentrated and self-alignment transmission. The architecture requires the UD to participate in channel formation, enabling dynamic adaptation to arbitrary movements of the UD while maintaining high efficiency.
    \item We establish a comprehensive theoretical model to analyze the proposed system based on EM propagation theory and a self-reproduction model, accounting for the noise introduced during the iterative process and establishing the corresponding transition matrices at each step. The model reveals the underlying mechanisms of resonance-driven transmission and demonstrates the system's ability to overcome challenges such as sidelobe interference and beam tracking limitations. 
    \item Numerical evaluations show that the proposed system achieves high-efficiency transmission and high-spectral-efficiency communication in the near field without additional tracking control. Furthermore, increasing the number of antennas at either the BS or the MT extends the effective operating range of RF-RBS without increasing the convergence time required for positioning and tracking the MT. 
\end{itemize}

The remainder of this article is organized as follows. Section II describes the RF-RBS architecture, RAAs, and resonance mechanism. Section III establishes an analytical model to calculate power flow, transmission efficiency, and communication performance. Section IV presents a numerical analysis to demonstrate the superior performance of the RF-RBS compared to the RD-BFS. Finally, Section V concludes the paper.

\begin{table}
    \caption{PRIMARY SYSTEM PARAMETERS}
    \centering
    {
        \begin{tabular}{>{\centering\arraybackslash}p{0.6in}>{\centering\arraybackslash}p{2.4in}}
        \hline
        \textbf{Symbol}                        & \textbf{Explanation}   
        \\ \hline
        $c$                         & Speed of EM wave  \\
        $f_c$, $\lambda$            & Carrier frequency and Wavelength \\
        $Z_0$                       & Characteristic impedance  \\
        $\alpha$, $\beta$           & Path-loss exponent and Scaling factor   \\
        
        $E$, $s$                   & Signal electric field and Baseband signal \\
        $p$, $v$, $\phi_0$          & Signal power, amplitude, and initial phase  \\

        $G_\text{t}$, $G_\text{r}$  & Antenna gain  \\
        $L$, $L_{nm}$               & Distance between transmitter and receiver  \\

        $f_\text{LO}$, $v_\text{LO}$    & LO signal frequency and amplitude  \\
        $G_\text{a}$, $\phi_\text{a}$   & PA gain and phase-lag  \\

        $k$                         & Iterations  \\
        $N$, $M$                    & BS and MT antenna elements\\

        $\mathbf{H}_\text{c}$       & Channel gain matrix \\
        $\mathbf{H}_\text{p}$       & Phase-conjugate transmission matrix \\
        $\mathbf{n}$, $F$           & Additive noise and Noise figure      \\

        $v_\text{m}$                & Limiter maximum signal amplitude     \\
       
        $\iota$, $\gamma$           & Iteration power loss and BS power gain  \\
        $P_\text{BS}$               & BS total radiated power  \\
        $P_\text{MT}$               & MT total received power  \\
        
        $W$                         & Bandwidth        \\
        $\alpha_{\text{pd}}$        & Feedback ratio     \\
        $\tilde{C}$                 & Spectral efficiency 
        \\ \hline
    \end{tabular}
    }
    \label{table: primary variable}
\end{table}

Notation: Boldface lower-case letters are vectors (e.g., $\mathbf{s}$), whereas boldface upper-case letters are matrices (e.g., $\mathbf{H}$). $\mathbf{I}_N$ is the identity matrix of size $N$, $h_{n,m} = [\mathbf{H}]_{n,m}$ represents the $(n, m)$-th element of matrix $\mathbf{H}$. The notation $\mathbf{x} \sim \mathcal{CN}(m, \sigma^2)$ indicates a complex circular symmetric Gaussian random variable (RV) with mean $m$ and variance $\sigma^2$, whereas $\mathbf{x} \sim \mathcal{CN}(m, \mathbf{C})$ denotes a complex Gaussian random vector with mean $m$ and covariance matrix $\mathbf{C}$. Constant $j = \sqrt{-1}$ denotes the imaginary unit. The real part of a complex number $z$ is $\Re(z)$.

Table~\ref{table: primary variable} includes the primary system parameters utilized in the subsequent analysis.

\section{System Overview}
\label{sec:sys}

RF-RBS is a promising technology that enables self-alignment and efficient RF band transmission based on RAAs and the resonance mechanism. In this section, we first introduce the system design of the RF-RBS. Then, we provide a detailed description of the RAA structure. Finally, the resonance mechanism is introduced.

\subsection{System Design of RF-RBS}

\begin{figure*}
\centering
\includegraphics[width=6.8in]{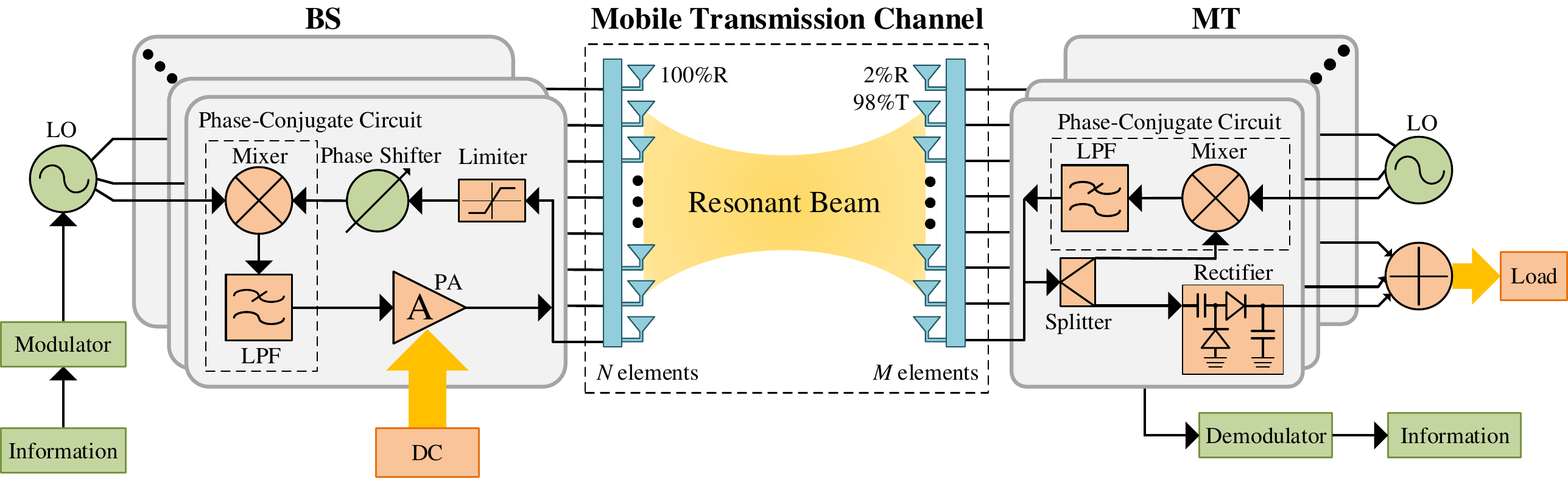}
\caption{Transceiver structure of RF-RBS.}
\label{fig: structure of RF-RBS}
\end{figure*}

Figure~\ref{fig: structure of RF-RBS} illustrates the transceiver design of the RF-RBS, including a BS and an MT operating in the near-field region, forming a mobile transmission channel.

The BS comprises an antenna array of $N$ elements with phase-conjugate circuits and power amplifiers (PAs). Each element in the BS includes an amplitude limiter to ensure the received power stays within a specific range. Following this, a phase shifter balances the phase difference between the BS and MT due to different circuit paths. Information to be transmitted is first modulated in the local oscillator (LO) and then mixed with the input signal in the phase-conjugate circuit. Finally, the phase-conjugated waves are amplified and re-emitted.

The MT includes an antenna array of $M$ elements with phase-conjugate circuits. A power divider directs most energy to subsequent energy and communication parts, while the remaining signal is returned to the antenna through a phase-conjugate circuit. The output signal is then rectified, summed, and delivered to the load in the energy processing unit. Information is extracted by demodulating the phase of the received signal from the central element. 

The mobile transmission channel is established between the BS and MT. Power and information flow in the RF-RBS occurs as follows:
1) DC Source: The DC source supplies power to the RF PAs in the BS to amplify EM waves.
2) Power and Information Transfer: The EM waves propagate back and forth in free space, transferring power and information from the BS to the MT. The MT sends a portion of the power back to the BS, where it is further amplified. 
3) Power and Information Output: The remaining power at the MT is converted to DC power via a rectifier circuit and supplied to a load, and the information is demodulated from the signal.

The system design ensures efficient power transfer and communication from the BS to the MT, utilizing resonance mechanisms to enhance system performance and reliability.

\subsection{Retro-Directive Antenna Array}

The RAA retransmits waves back to their original direction using phase conjugation as the fundamental principle. Phase conjugation can be achieved through both digital and analog methods. The digital approach involves detecting and calculating the phase of the incident wave before emitting a phase-conjugate wave. Conversely, the analog method offers a faster and more cost-effective means of achieving phase conjugation. For RF-RBS, adopting a heterodyne technique for phase conjugation, as demonstrated by the Pon Array \cite{PonRetrodirective1964}, is an effective strategy to enable retro-directivity. Other techniques, including metasurfaces \cite{WongBinary2018, abbasi2019maxwell}, have also shown potential effectiveness.

Figure~\ref{fig: structure of RF-RBS} presents a typical design of the RAA at both the BS and MT, based on the Pon Array. The antenna initially processes the incoming EM wave using a mixer that interacts with a LO of frequency $f_\text{LO} = 2 f_c$. The resulting mixed wave is filtered to produce a phase-conjugate EM wave at the frequency $f_c$. Specifically, at the BS, the phase-conjugate wave is subsequently amplified by a power amplifier (PA) and radiated from the same antenna. At the MT, a portion of the incoming power is mixed and filtered in the phase-conjugate circuit to generate a phase-conjugate EM wave emitted from the same antenna.

Integrating the transmitting and receiving antennas is crucial to enable a single antenna to function in both modes. In digitally implemented phase-conjugate antennas, a switching module toggles the antenna's operating mode. A time-division strategy controls whether it operates in transmitting or receiving modes. Alternatively, in the analog approach, phase-conjugate antennas can be arranged into an RAA \cite{bird2005design}. This configuration utilizes dual-polarized antennas, where the input and output EM waves exhibit different polarizations \cite{zhou2014retrodirective, kang2023simultaneous}, allowing the antenna to operate simultaneously in transmitting and receiving modes.

\subsection{Resonance Mechanism}

\begin{figure}
\centering
\includegraphics[width=3.4in]{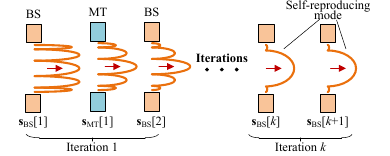}
\caption{Round-trip iteration process and self-reproducing mode in RF-RBS.}
\label{fig: Fig-RBS-Iteration}
\end{figure}

Based on the RAA, the RD-BFS is typically used for directional transmission of EM waves from the BS to the MT. 
The RD-BFS utilizes the current pilot signal from the MT to direct the transmission of EM waves from the BS to the MT. Each element in the BS emits phase-conjugated and power-amplified EM waves, focusing energy on the MT. However, this method ensures the highest received power only for the central antenna on the MT, where the pilot signal originates, resulting in unpredictable received power for other antennas. This situation can lead to increased sidelobes and diminished efficiency, as noted in \cite{chou2020conformal}.

In contrast, the RF-RBS achieves efficient EM wave transmission between the BS and MT through a resonance mechanism during the round-trip iteration process. As depicted in Fig.~\ref{fig: Fig-RBS-Iteration}, after multiple round trips of EM waves between the BS and MT, the field distribution stabilizes and becomes self-reproducing, as the resonant field forms through the superposition of in-phase waves, canceling out other phases. By achieving resonance, the RF-RBS maximizes the efficiency of EM wave transmission and maintains high-capacity communication links between the BS and MT. In detail, resonance in RF-RBS allows the system components to operate at the optimal frequency. At the same time, it enhances energy control, ensuring that EM wave distribution is optimized between the BS and MT, achieving efficient transmission. Finally, the resonance phenomenon enhances signal transmission intensity and inhibits interference.

To achieve resonance in RF-RBS, several conditions must be met: appropriate frequency, necessary boundary conditions, and a matched excitation source. The system's operating frequency should be within a specific resonant frequency range, with all components tuned to a uniform single frequency $f_c$. The boundary conditions must allow for proper reflection of EM waves, enabling multiple reflections to form a stable standing wave mode. In RF-RBS, RAAs enable EM wave transmission between BS and MT and ensure phase alignment through phase shifters. Finally, the excitation source must provide sufficient resonance power. In RF-RBS, the PAs supply the required energy, compensating for transmission, device, and output losses. Other devices, such as masers \cite{OxborrowRoommaser2012}, which use stimulated radiation and resonance mechanisms for microwave amplification, are also viable options.

The single-frequency design is specifically proposed in RF-RBS. Compared to the dual-frequency design \cite{Kang2023Dualfrequency}, it eliminates the need for additional arrays with different element spacing and ensures the collinearity of the round-trip beams, further enhancing beam alignment accuracy and improving EM safety. However, a challenge arises from the phase delay difference between the BS and MT, which prevents the formation of standing waves in the mobile channel. As mentioned earlier, phase shifters in the BS are used to achieve the required phase alignment. Additionally, the design and manufacturing of high-gain, stable power amplifiers and phase shifters require further research.

Owing to the antenna array's propagation characteristics, uncorrelated noise is effectively attenuated during the beamforming process. Additionally, the resonance mechanism helps eliminate out-of-phase waves introduced by noise \cite{aokiObservationStrongCoupling2006}, keeping noise power low. The established resonance ensures that the system maintains a high signal-to-noise ratio (SNR).

\section{Analysis Model}
\label{sec:ana}

In this section, we develop an analytical model to explore the dynamics of EM wave flow and propagation within the RF-RBS and calculate transmission efficiency, resonance establishment, and communication capacity.

\subsection{EM Waves Radiation and Propagation}

\subsubsection*{Power Radiation}

The electric field component of an EM wave at a given time $t$ is described by
\begin{equation}
    E(t) = \Re\{s(t) e^{ j 2 \pi f_c t}\},
\end{equation}
where $s(t)$ represents the baseband signal, $f_c$ is the carrier frequency of the signal. The averaged effective power radiated by an antenna is then given by
\begin{equation}
    p(t) = \frac{|s(t)|^2}{2Z_0},
\end{equation}
where $Z_0$ is the characteristic impedance of the medium. Therefore, the radiated signal from an antenna can be modeled as
\begin{equation}
    s(t) = \sqrt{2 Z_0 p(t)} \cdot e^{j \phi_0(t)},
\end{equation}
where $p(t)$ and $\phi_0(t)$ are the averaged effective radiated power and the initial radiated phase of the antenna, respectively. $v(t) = \sqrt{2 Z_0 p(t)}$ is the amplitude of the signal.

\subsubsection*{Power Propagation}

A basic wireless channel model between the power transmitter and receiver includes large-scale path loss and small-scale multipath fading. The channel gain for an RF signal is expressed as \cite{khang2018microwave, choi2018distributed}
\begin{equation}
    h_\text{c} = \sqrt{\beta} L^{-\frac{\alpha}{2}} \sqrt{G_\text{t} G_\text{r}} e^{j \frac{2 \pi}{\lambda} L},
\end{equation}
where $\lambda = c/f_c$ is the wavelength and $c$ is the speed of the EM wave in free space. $\beta > 0$ is a constant scaling factor, $L$ is the distance between the transmitter and receiver, $\alpha \geq 2$ denotes the path-loss exponent, and $G_{\text{t}}$, $G_{\text{r}}$ are the gains of the transmitting and receiving antennas, respectively. The antenna's gain is related to the characteristics of the antenna and changes with the change of angle. Using the Friis transmission equation, we have $\beta = \frac{\lambda^2}{16\pi^2}$ and $\alpha = 2$.

To model multi-antenna signal transmission from $N$ transmitting antennas to $M$ receiving antennas in the near-field region, the transmission channel gain matrix $\mathbf{H}_\text{c} \in \mathbb{C}^{M \times N}$ is introduced, and for each element $h_{\text{c}, n, m}$ in $\mathbf{H}_\text{c}$ is
\begin{equation}
    h_{\text{c}, n, m} = \sqrt{\beta} L_{nm}^{-\frac{\alpha}{2}} \sqrt{G_{\text{t},n} G_{\text{r},m}} e^{j \frac{2 \pi}{\lambda} L_{nm}},
\end{equation}
where $L_{nm}$ denotes the path length from the $n$-th transmitting antenna to the $m$-th receiving antenna. $G_{\text{t},n}$ and $G_{\text{r},m}$ are the gains of the transmitting and receiving antennas, respectively. As the BS and MT operate in a near-field environment, EM wave propagation adheres to spherical wavefronts, requiring the phase differences across the antenna arrays to be rigorously calculated and accurately derived.

Consequently, the received signal of an antenna array is modeled as
\begin{equation}
    \mathbf{s}_\text{r}(t) = \mathbf{H}_\text{c} \mathbf{s}_\text{t}(t) + \mathbf{n}_\text{c},
    \label{eq: power transmission}
\end{equation}
where $\mathbf{s}_\text{t} \in \mathbb{C}^{N \times 1}$ and $\mathbf{s}_\text{r} \in \mathbb{C}^{M \times 1}$ are the transmitted and received signal vectors, respectively. $\mathbf{n}_\text{c} \in \mathbb{C}^{M \times 1}$ is the additive white Gaussian noise (AWGN), with $\mathbf{n}_\text{c} \sim \mathcal{CN}(0, \sigma_{\text{c}}^2 \mathbf{I}_M)$. Note that $\sigma_{\text{c}}^2 = 2Z_0 \kappa T_0 F_\text{c} W$, being $\kappa$ the Boltzmann constant, $T_0 = 290 \, \text{K}$, $F_\text{c}$ the noise figure during receiving signal, and $W$ is the bandwidth of the narrow signal \cite{dardari2023establishing}. No antenna mutual coupling is assumed to occur between antennas because at least half a wavelength separates all antennas.  

\subsection{Retro-Directive Antenna Array}

An RAA retransmits the EM wave back to its source using phase conjugation as the primary mechanism. This process, illustrated in Fig.~\ref{fig: structure of RF-RBS}, involves mixing the input signal, $E_\text{pi}(t) = v_\text{pi} \cos(2 \pi f_c t + \phi_\text{pi})$, with a LO signal, $E_\text{LO}(t) = v_\text{LO} \cos(2 \pi f_\text{LO} t + \phi_\text{LO})$, where $v_\text{pi}$ and $v_\text{LO}$ are the amplitudes, $f_\text{LO} = 2 f_c$ is the frequency of LO signal, and $\phi_\text{pi}$ and $\phi_\text{LO}$ are the phases of the input and LO signals, respectively. The resulting mixed signal is
\begin{equation}
    \begin{aligned}
        E_\text{mix}(t) = & \frac{v_\text{pi} v_\text{LO}}{2} \left[\cos(2 \pi (f_\text{LO} - f_c) t - \phi_\text{pi} + \phi_\text{LO} ) \right. \\
            & \left. + \cos \left( 2 \pi (f_\text{LO} + f_c) t + \phi_\text{pi} + \phi_\text{LO} \right) \right].
    \end{aligned}
\end{equation}
Filtering out the high-frequency component $(f_\text{LO} + f_c)$ through the LPF, the intermediate frequency (IF) component is retained as
\begin{equation}
    E_\text{po}(t) = \frac{v_\text{pi} v_\text{LO}}{2} \cos(2 \pi f_c t - \phi_\text{pi} + \phi_\text{LO}),
\end{equation}

Therefore, the phase of the output signal after the phase-conjugate circuit is effectively the conjugate of the input phase
\begin{equation}
    \phi_\text{po} = -\phi_\text{pi} + \phi_\text{LO},
\end{equation}
where $\phi_\text{LO}$ includes any constant phase delay introduced by the phase-conjugate circuit. Additionally, the phase-conjugate circuit introduces amplitude attenuation of $\frac{v_\text{LO}}{2}$.

Using the transmission matrix $\mathbf{H}_\text{p} \in \mathbb{C}^{N \times N}$, the relationship between the output and input signals for a phase-conjugate array with $N$ elements is therefore
\begin{equation}
    \mathbf{s}_\text{po}(t) = \mathbf{H}_\text{p} \mathbf{s}_\text{pi}(t) + \mathbf{n}_\text{p},
\end{equation}
where $\mathbf{s}_\text{pi}(t)  \in \mathbb{C}^{N \times 1}$ and $\mathbf{s}_\text{po}(t)  \in \mathbb{C}^{N \times 1}$ are input and output signal vectors, respectively. $\mathbf{n}_\text{p} \in \mathbb{C}^{N \times 1}$ is the AWGN, with $\mathbf{n}_\text{p} \sim \mathcal{CN}(0, \sigma_{\text{p}}^2 \mathbf{I}_N)$. Note that $\sigma_{\text{p}}^2 = 2Z_0 \kappa T_0 F_\text{p} W$, being $F_\text{p}$ the noise figure in the phase-conjugate circuit \cite{zhang2022active}. The transmission matrix $\mathbf{H}_\text{p}$ of the phase conjugate circuit is defined as $\text{diag}(h_{\text{p},1}, h_{\text{p},2}, \ldots, h_{\text{p},N})$, and each element is
\begin{equation}
    h_{\text{p},n} = \frac{v_\text{LO}}{2} e^{-j2\phi_{\text{pi}, n} + 
     \phi_\text{LO}},
\end{equation}
where $\phi_{\text{pi}, n}$ is the phase of the $n$-th element prior to conjugation. 

\subsection{Power Flow in RF-RBS}

In RF-RBS, power flow is bidirectional due to the resonance mechanism, allowing some power to circulate between the BS and MT. Defining one complete cycle of power—from the BS to the MT and back—as one iteration (denoted by $k$), and ignoring the superposition of noise, the EM wave signal received by the MT for the $(k+1)$-th iteration is given by
\begin{equation}
    \mathbf{s}_\text{MT}[k+1] = \mathbf{H}_\text{cT} \mathbf{H}_{\text{BS}} \mathbf{H}_\text{cR} \mathbf{H}_{\text{MT}} \mathbf{s}_\text{MT}[k],
    \label{eq: signal circulation}
\end{equation}
where $\mathbf{s}_\text{MT}[k]$ is the received signal at the MT during the $k$-th iteration. $\mathbf{H}_\text{cT} \in \mathbb{C}^{M \times N}$ and $\mathbf{H}_\text{cR} \in \mathbb{C}^{N \times M}$ are the transmission channel gain matrices for EM waves emitted from the BS and MT, respectively. $\mathbf{H}_{\text{BS}}  \in \mathbb{C}^{N \times N}$ and $\mathbf{H}_{\text{MT}} \in \mathbb{C}^{M \times M}$ are the transmission matrices of the circuits in the BS and MT, respectively.

At the BS, the received signal undergoes several processing stages before being emitted by the antenna. Let matrices $\mathbf{H}_\text{p} \in \mathbb{C}^{N \times N}$ and $\mathbf{H}_\text{a} \in \mathbb{C}^{N \times N}$ represent the transmission matrices of the phase-conjugate circuits and PAs, respectively. These matrices are diagonal, indicating they only act on their respective antenna elements. The cumulative transmission matrix $\mathbf{H}_\text{BS}$ at the BS is given by
\begin{equation}
    \mathbf{H}_\text{BS} = \mathbf{H}_\text{a} \mathbf{H}_\text{p} e^{j \phi_\text{s}} \alpha_{\text{l}}.
\end{equation}
where $\alpha_{\text{l}}$ is the attenuation ratio of the limiter, ensuring that the amplitude of the output signal is less than the allowed maximum signal amplitude $v_\text{m}$, allowing subsequent circuits to operate normally without exceeding the power limit. The $\phi_\text{s}$ is the phase-lag value of the phase shifter, a fixed value set according to the circuit characteristics.

The PAs produce the RF power required by the components and power output in RF-RBS. The bandwidths of PAs are typically much larger than the working bandwidth, so we neglect the dynamics of PAs. Each element in $\mathbf{H}_\text{a}$ can be modeled as \cite{moghadam2018energy}
\begin{equation}
    h_{\text{a},n} = \sqrt{G_\text{a}} \cdot e^{j \phi_\text{a}},
\end{equation}
where $G_\text{a}$ is the power gain, and $\phi_\text{a}$ is the phase-lag in the PA. An RF PA may be non-linear, with its gain as a function of the input signal level, often denoted the amplitude-to-amplitude (AM/AM) distortion. Generally, PAs behave linearly only within a particular input power range, with nonlinearity becoming significant at larger input power levels. Inevitably, PAs introduce additional noise. Let $n_\text{a}$ be the AWGN for each PA, with $n_\text{a} \sim \mathcal{CN}(0, \sigma_{\text{a}}^2)$. Note that $\sigma_{\text{a}}^2 = 2Z_0 \kappa T_0 F_\text{a} W$, where $F_\text{a}$ is the noise figure of the PA.

At the MT, the signals pass through the power dividers and phase-conjugate circuits. The transmission matrix $\mathbf{H}_\text{MT}$ at the MT is given by
\begin{equation}
    \mathbf{H}_\text{MT} = \mathbf{H}_\text{p} \alpha_{\text{pd}} ,
\end{equation}
where $\alpha_{\text{pd}}$ is the ratio of the power divider entering the subsequent phase-conjugate circuit. The remaining $1-\alpha_{\text{pd}}$ goes into subsequent information and energy processing units.

\subsection{Resonance Establishment}

\begin{figure}
    \centering
    \includegraphics[width=3.4in]{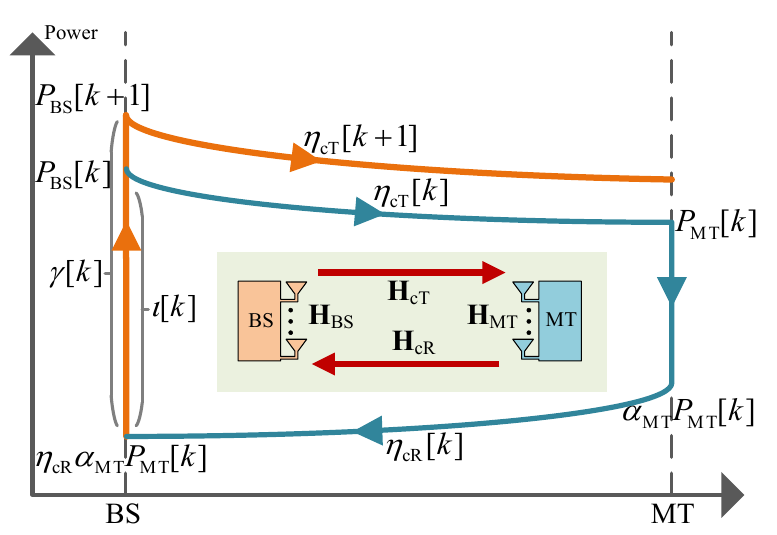}
    \caption{Power circulation between BS and MT.}
    \label{fig: power circulation}
\end{figure}

During the power cycle, losses occur due to transmission, reception, and circuit processing, which result in decreased power within the system. PAs at the BS amplify the input signal, thereby increasing the signal power. The power circulation for a specific iteration $k$ between the BS and the MT is illustrated in Fig.~\ref{fig: power circulation}. In iteration $k$, the BS transmits a signal with a total power of $P_\text{BS}[k] = \sum_{n=1}^N {p}^\text{BS}_{\text{t},n}[k]$. The power received on the MT is $P_\text{MT}[k] = \sum_{m=1}^M {p}^\text{MT}_{\text{r},m}[k]$, with the transmission efficiency $\eta_\text{cT}[k] = P_\text{MT}[k]/P_\text{BS}[k]$ from the BS to MT.

After reception, the MT returns a small amount of power to the BS, scaled by the factor $\alpha_\text{MT}$ through a power splitter and phase-conjugate circuit. Thus, in iteration $k+1$, the transmission power of the BS is
\begin{equation}
    P_\text{BS}[k+1] = G_\text{BS}\left(\eta_\text{cR} \alpha_\text{MT} P_\text{MT}[k]\right),
\end{equation}
where $G_\text{BS}(\cdot)$ represents the total gain function of the BS, including $G_\text{a}$ of PAs and other circuits. In the iteration $k$, the difference between the power emitted from the BS and that subsequently received is denoted as $\iota[k]$, and the power gain as $\gamma[k]$, defined mathematically as
\begin{equation}
    \begin{cases}
    \iota[k] = P_\text{BS}[k] - \eta_\text{cR} \alpha_\text{MT} P_\text{MT}[k] = (1 - \alpha_\text{MT} \eta_\text{cR} \eta_\text{cT}) P_\text{BS}[k] \\
    \gamma[k] = P_\text{BS}[k+1] - \alpha_\text{MT} \eta_\text{cR} \eta_\text{cT} P_\text{BS}[k]
    \end{cases} .
    \label{eq: Gain loss definition}
\end{equation}

When resonance is established, $\iota$ and $\gamma$ should converge, i.e., the power lost in transmission and harvested by the MT is exactly compensated by the PAs. If resonance has not been achieved, the power gain should always exceed the loss, prompting the system to transmit with increased power to reach a stable state. The process is formally described as  
\begin{equation}
    \begin{cases}
        \gamma[k=1] > \iota[k=1] \\
        \gamma[k \to \infty] = \iota[k \to \infty]
    \end{cases} .
    \label{eq: power ratio selection}
\end{equation}

The system eventually achieves stability when the above gain and loss conditions are met. The RF-RBS enters a self-reproducing mode when the field distribution between the BS and MT arrays stabilizes, signifying constant transmission efficiency and power. Assuming the system converges at iteration $k$, according to the resonance principle, the field distribution in iteration $k + 1$ remains unchanged, which can be described as 
\begin{equation}
    \mathbf{s}_\text{MT}[k+1] = \mathbf{s}_\text{MT}[k].
    \label{eq: resonance state}
\end{equation}
The effect of noise on the signal is ignored here. In fact, due to the addition of noise in the above iterative process, the signal will still have a slightly unpredictable change, but the whole remains unchanged. Upon validating \eqref{eq: resonance state}, it is confirmed that stable resonance is achieved, and the fields reproduce themselves in subsequent iterations. 

\subsection{Energy Harvesting Process}

\begin{figure}
    \centering
    \includegraphics[width=3.3in]{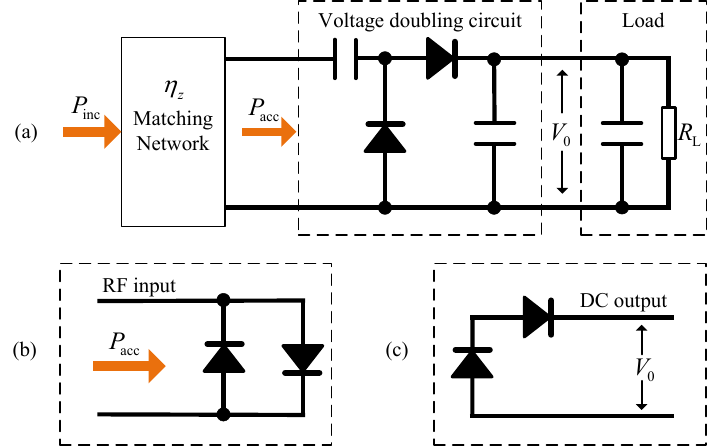}
    \caption{Single energy harvesting circuit. (a) Voltage doubling circuit. (b) RF equivalent circuit. (c) DC equivalent circuit.}
    \label{fig: rectenna}
\end{figure}

In each element of the MT, most of the received RF power is directed to energy harvesting circuits, which provide DC power for the load. To facilitate RF-to-DC conversion, each element is outfitted with a rectifier and is integrated using a combiner, as depicted in Fig.~\ref{fig: structure of RF-RBS}. Before rectifier, a matching network ensures effective impedance matching, assumed to have achieved an impedance matching efficiency of $\eta_\text{z}$.

The voltage doubling circuit is employed to enhance the DC output voltage and conversion efficiency, as illustrated in Fig.~\ref{fig: rectenna}(a). When an RF signal is introduced to the circuit, capacitors act as short circuits, leading to the equivalent circuit shown in Fig.~\ref{fig: rectenna}(b). This configuration effectively halves the RF input impedance relative to a single diode. For DC signals, capacitors behave as open circuits, allowing the two diodes to function as two voltage sources connected in series, as shown in Fig.~\ref{fig: rectenna}(c).

The implicit expression for a diode's output voltage $V_0$ is given by \cite{visser2013RFEnergyHarvesting}
\begin{equation}
    I_0 \left( \frac{q}{n_0k_0T} \sqrt{8 R_g P_{\text{acc}}} \right) = \left( 1 + \frac{V_0}{R_\text{L} I_\text{s}} \right) e^{ \left( 1 + \frac{R_\text{g} + R_\text{s}}{R_\text{L}} \right) \frac{q}{n_0k_0T} V_0 }
    \label{eq: diode expression}
\end{equation}
where $q$ represents the electron charge, $n_0$ is the diode’s ideality factor, $k_0$ denotes Boltzmann’s constant, and $T$ is the temperature in Kelvin. $I_0(\cdot)$ is the zero-order modified Bessel function of the first kind, and $I_\text{s}$ is the diode’s saturation current. $R_\text{L}$, $R_\text{g}$, and $R_\text{s}$ denote the load resistance, input resistance, and diode series resistance, respectively. $P_{\text{acc}} = \eta_\text{z} P_{\text{inc}}$ represents the power accepted after matching, with $P_{\text{inc}}$ being the incident RF power. 

To apply \eqref{eq: diode expression} of a single diode to the voltage doubling circuit in Fig. \ref{fig: rectenna}, $R_\text{g}$ is doubled, and $R_\text{L}$ is halved to maintain the same current inflow and outflow from the circuit \cite{Agilent1999}. Once the output voltage $V_0$ is determined, the conversion efficiency of the voltage doubling circuit is calculated as
\begin{equation}
    \eta_\text{con} = \frac{V_0^2}{R_\text{L} P_{\text{acc}}}.
\end{equation}

After converting the RF power to DC power in each element, a DC combiner sums all the DC outputs, excluding the centrally selected element used for communication. Unlike an RF combiner, a DC combiner does not require a complex phase control circuit. Thus, the total output power delivered to the load is derived by summing these outputs
\begin{equation}
    P_{\text{dc}} = \sum_{\substack{m=1 \\ m \neq m_c}}^{M} \eta_\text{con} \eta_\text{z} (1-\alpha_{\text{pd}}) {p}^\text{MT}_{\text{r},m}
\end{equation}
where $m_c$ signifies the antenna selected for information demodulation, and ${p}^\text{MT}_{\text{r},m}$ is the received power of the $m$-th element in MT.

\subsection{Communication Channel}

In RF-RBS, downlink communication begins only after stable resonance is established, ensuring an optimized transmission path with minimal interference and maximum signal strength. The BS uses a phase modulator to encode information onto the LO signal and transmit it along the resonant path. Various phase-based signaling schemes, such as Quadrature Phase Shift Keying (QPSK), can be utilized in RF-RBS. The received signal from the central antenna is selected at the MT, and the information is extracted through subsequent demodulation. Additionally, the communication channel in RF-RBS can be regarded as a linear time-invariant system \cite{xiong2021retro}.

The received signal at the MT contains noise generated during iterations. However, because stable resonance is established in RF-RBS, the transmitted power is generally much greater than the thermal noise, and the antenna array's self-noise improvement effect further reduces the noise. Therefore, in the analysis, the noise introduced by previous iterations is ignored, and only the noise of the current iteration is considered in communication.

Considering the current iteration, at the BS, the received signal from the MT is combined with receiving noise, which is further combined with the noise from the phase-conjugate circuit. In the PAs, the noise and the signal are amplified simultaneously, and the PA noise is added. Therefore, the noise variance of each element in the BS is
\begin{equation}
    \sigma_{\text{BS}, n}^2 =  G_\text{a} \sigma_{\text{c}}^2 + G_\text{a} \sigma_{\text{p}}^2 + \sigma_{\text{a}}^2 ,
\end{equation}

After the MT receives the signal, it passes through a power divider, which directs a proportion $(1 - \alpha_{\text{pd}})$ of the power to subsequent information and energy processing units. The received signal power of the selected central antenna is ${p}^\text{MT}_{\text{r},m_c}$, and the SNR of the system is given by
\begin{equation}
    \text{SNR} = \frac{(1 - \alpha_{\text{pd}}) 2 Z_0 {p}^\text{MT}_{\text{r},m_c}}{(1 - \alpha_{\text{pd}}) \boldsymbol{\sigma}_{\text{BS}}^2 + \sigma_{\text{c}}^2 + 2 Z_0 \kappa T_0 F_\text{d} W} ,
\end{equation}
where $F_\text{d}$ is the noise figure of the demodulation circuits.

The spectral efficiency of the RF-RBS for downlink communication can be computed according to Shannon’s theory \cite{Akdeniz2014Millimeter, LTE2007Mogensen}
\begin{equation}
    \tilde{C} = \log_2 \left\{ 1 + 10^{0.1(\text{SNR} - \Delta)} \right\} ,
\end{equation}
with SNR and channel loss factor $\Delta$ given in dB.

\section{Performance Evaluation}
\label{sec:perfrom}

\begin{table}
    \caption{PARAMETERS IN THE SIMULATION SETTING}
    \label{table: system parameters}
    \centering
    \begin{tabular}{ccc}
    \hline
        \textbf{Symbol}         & \textbf{Parameter}            & \textbf{Value}    \\
    \hline
        $f_c$                   & Carrier frequency             & 30~GHz            \\
        $\lambda$               & Wavelength                    & 1~cm              \\
        $d$                     & Spacing between adjacent antennas & $\frac{\lambda}{2}$            \\
        $Z_0$                   & Characteristic impedance      & 50 $\Omega$       \\ 
        $\alpha$                & Path-loss exponent            & 2                 \\
        $\beta$                 & Scaling factor                & $\frac{\lambda^2}{16\pi^2}$ \\
        $W$                     & Bandwidth                     & 500~MHz                 \\

        $G_\text{t}, G_\text{r}$& Antenna gain & 5 dBi (Max) \cite{balanis2016antenna}   \\
        $G_\text{a}$            & PA gain & 20 dB (Max)   \\
        $\phi_\text{a}$         & PA phase-lag & $\frac{\pi}{6}$       \\
        
        $N$                     & BS antenna elements       & $40 \times 40$, $50 \times 50$          \\
        $M$                     & MT antenna elements       & $40 \times 40$, $50 \times 50$       \\
        $v_\text{LO}$           & LO signal amplitude              & 2~V   \\

        $v_\text{m}$            & Limiter maximum signal amplitude     & $\sqrt{0.2}$~V \\

        $F_\text{c}$            & Antenna noise figure      & 3~dB \cite{dardari2023establishing}     \\
        $F_\text{p}$            & Phase-conjugate circuit noise figure      & 6~dB \cite{dardari2023establishing}     \\
        $F_\text{a}$            & PA noise figure           & 5~dB \cite{analog_hmc1132}     \\
        $F_\text{d}$            & Demodulation circuits noise figure           & 7~dB \cite{Akdeniz2014Millimeter}     \\

        $\alpha_{\text{pd}}$    & Power feedback ratio              & 0.02   \\
        $\Delta$        & Channel loss              & 3~dB \cite{A2020Dutta}   \\

        $\eta_\text{z}$      & Matching efficiency     & 0.95   \\
        $R_\text{L}$        & Load resistance          & 100 $\Omega$   \\
        $R_\text{g}$       & Input resistance          & 50 $\Omega$   \\
        $R_\text{s}$        & Diode series resistance             & 25 $\Omega$   \\
        $I_\text{s}$        & Diode’s saturation current             & $10^{-6}$ A   \\
        
    \hline
    \end{tabular}
\end{table}

This section analyzes the dynamic optimization process, transmission efficiency, spatial power distribution, and communication performances of the proposed RF-RBS based on the analytical model. Additionally, we compare the performances of RF-RBS and RD-BFS under varying distances and horizontal offsets.

\subsection{Parameter Setting}

\begin{figure}
    \centering
    \subfloat[$\phi=0^\circ/180^\circ$]{%
    \includegraphics[width=1.6in]{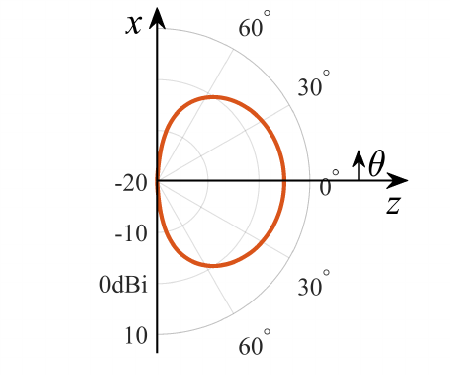}
    } 
    \subfloat[$\phi=90^\circ/270^\circ$]{%
        \includegraphics[width=1.6in]{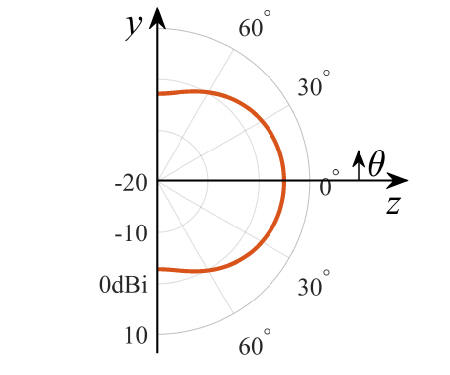}
    }
    \caption{Antenna gain versus elevation angle $\theta$ and azimuth angle $\phi$.}
    \label{fig: Antenna Gain}
\end{figure}

\begin{figure}
    \centering
    \includegraphics[width=3.5in]{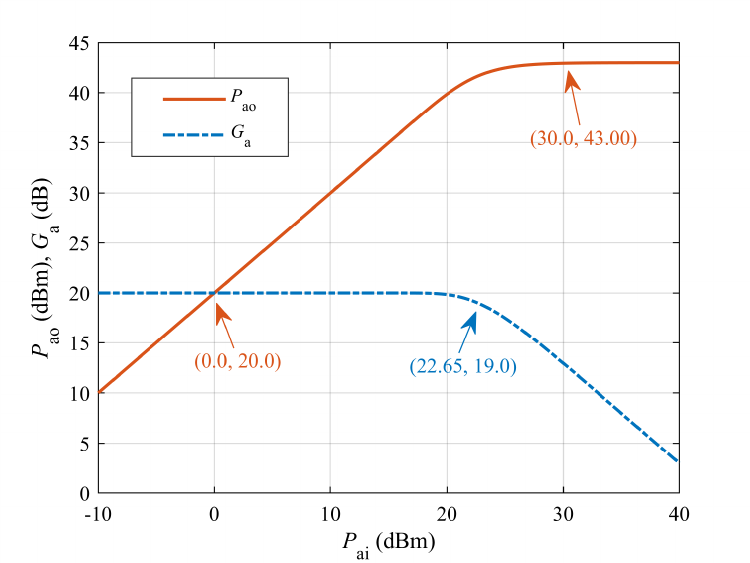}
    \caption{Output power and gain versus input power of RF PA.}
    \label{fig: RF PA}
\end{figure}

The parameters in the simulation setting are listed in Table~\ref{table: system parameters}. The RF-RBS operates at a frequency of $f_c = 30$~GHz, corresponding to a wavelength $\lambda$ of approximately 1~cm. Although 30~GHz is highlighted in this analysis, it is not the only feasible frequency. With appropriate antenna and component designs, a wide range of frequencies can be explored. To ensure optimal system performance, all subsystems—including the microstrip antennas, RF PAs, filters, and phase-conjugate circuits—should be meticulously tuned to work at 30~GHz.

The BS and MT have the same number of elements, with a typical spacing distance $d = \lambda/2$ between adjacent elements. For the antenna gain $G_{\text{t}}$ and $G_{\text{r}}$ of microstrip antennas, we consult the radiation pattern presented in formula (14-44) in \cite{balanis2016antenna}, which leads us to the antenna gain pattern shown in Fig.~\ref{fig: Antenna Gain}. Notably, the antenna gain attains its maximum when $\theta = 0$, equating to a value of $\pi$ or approximately 5~dBi. As $\theta$ increases, the antenna gain experiences a noticeable reduction. 

Figure~\ref{fig: RF PA} illustrates the gain characteristic of the PAs operating at 30~GHz based on the gain mathematical model \cite{Bogya2004Linear}. Within its standard operating range, the total output power $P_\text{ao}$ is exclusively determined by the total input power $P_\text{ai}$, obviating the need for additional controls. As depicted in the figure, the amplification gain $G_\text{a}$ initially remains high at 20~dB. It then gradually decreases as the input power increases until the output power stabilizes at approximately 43~dBm or 20~W.

In determining the power divider's feedback ratio $\alpha_{\text{pd}}$ at the MT, we adhere to the criterion outlined in equation \ref{eq: power ratio selection}. Based on the characteristics of the PA, we have set the feedback ratio as $\alpha_{\text{pd}} = 0.02$ to ensure optimal system operation. 

The following analysis positions the BS at the origin $(0, 0, 0)$ and the MT at $(\Delta x, 0, \Delta z)$, where $\Delta z$ indicates the vertical distance along the $z$-axis and $\Delta x$ represents the horizontal offset along the $x$-axis. The elevation angle $\theta$ is measured from the $z$-axis, with $\theta = 0$ along the $z$-axis. The azimuth angle $\phi$ is measured from the $x$-axis, with $\phi = 0$ along the $x$-axis. The elements at the BS and MT follow a rectangular grid distribution, with a spacing distance of $d$.

In the simulation, RF-BFS operates based on the model described above. Noise impacts the system’s convergence and stability, which is also introduced as the initial perturbation signal to start the iteration process in RF-BFS. The RD-BFS simulation below uses a setup identical in size and antenna gain to the RF-RBS configuration. Instead of iterating, a central antenna of the MT is selected to emit an isotropic pilot signal. The BS is configured to record the power and phase of the received signals and then emit a phase-conjugate response \cite{kim2023curved}. Since the radiated power of RD-BFS is configurable and lacks a comparable baseline, it is not displayed in the simulation results.

\subsection{Iteration Characteristics of RF-RBS}

\begin{figure}
    \centering
    \includegraphics[width=3.5in]{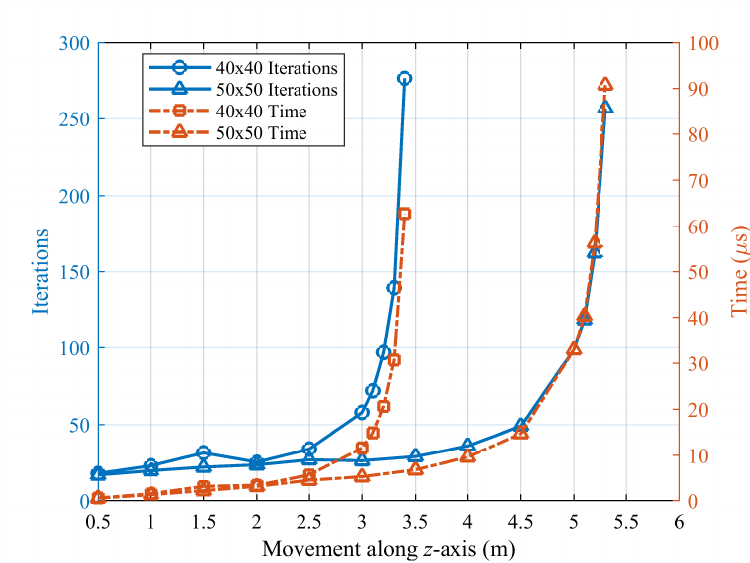}
    \caption{Iterations and time for system convergence versus MT's movement along z-axis $\Delta z$.}
    \label{fig: Fig-Iterations}
\end{figure}

\begin{figure}
    \centering
    \includegraphics[width=3.5in]{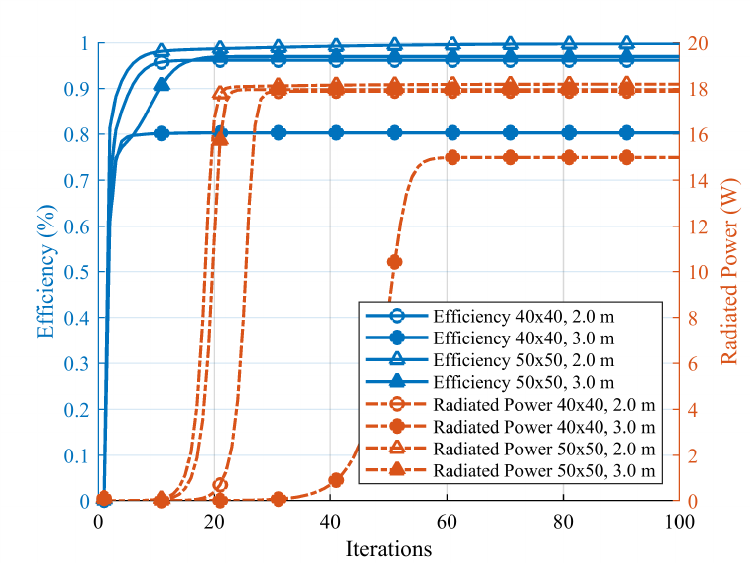}
    \caption{Spatial transmission efficiency $\eta_\text{c}$ and radiated power $P_\text{BS}$ versus the number of iterations in different distances.}
    \label{fig: Fig-Iteration-dynamic}
\end{figure}

\begin{figure*}
    \centering
        \subfloat[Spatial normalized power distribution at iterations 1, 2, 5, 10, and 100.]{%
        \includegraphics[width=7in]{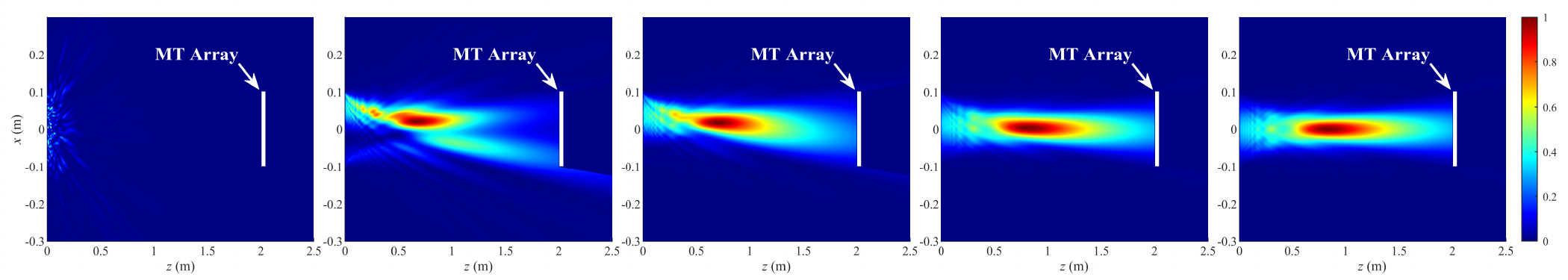}
        } \\
        \vspace{-3mm}
        \subfloat[Normalized power distribution on MT array at iterations 1, 2, 5, 10, and 100.]{%
        \includegraphics[width=7in]{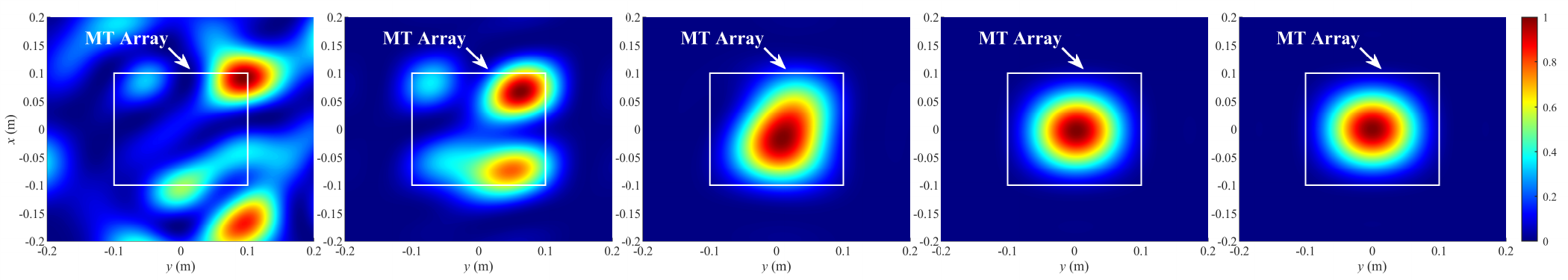}
        }
    \caption{Normalized power distribution on spatial $xOz$ plane and MT array across iterations.}
    \label{fig: Fig_spatial_dis}
\end{figure*}

RF-RBS completes beam optimization through continuous iteration, eventually realizing the self-reproducing mode of EM waves to form a resonant beam between the BS and MT. When the MT is within the effective working range of the BS, the EM waves will spontaneously transmit back and forth in the presence of noise, reaching a stable state.

Figure~\ref{fig: Fig-Iterations} shows the iterations and time required for the RF-RBS to converge at different distances and array sizes. In numerical simulations, system convergence is defined as when the power difference between two successive iterations is less than a specific threshold (set to 0.1\%). Considering the randomness of the noise, the data is averaged over 20 independent experiments. The figure shows that the $40 \times 40$ array only has data within 3.4~m, and the $50 \times 50$ array only has data within 5.3~m. This is because the RF-RBS cannot establish a stable EM beam beyond these distances, failing to meet the gain and loss conditions. Within their respective effective ranges, the number of iterations required increases as the distance increases and sharply rises as the maximum working distance approaches. Additionally, longer distances increase propagation time, increasing the time required for convergence. The effect of array size on the required iterations is insignificant because the optimization of RF-RBS is parallel for all elements.

Figure~\ref{fig: Fig-Iteration-dynamic} shows the changes in transmission efficiency $\eta_\text{cT}$ and radiated power $P_\text{BS}$ with iterations for different array sizes and distances. Initially, the power and efficiency are close to zero, as no array optimization and power amplification are carried out. As iterations progress, efficiency improves significantly while power remains low. This is due to the initial minimal power, where the PA gain is insufficient to produce a substantial increase. However, after a sufficient number of iterations, a noticeable rise in power becomes evident. Eventually, efficiency continues to increase to a stable level, and power increases rapidly at a certain point, reaching a stable level. It is shown that the longer distance lowers the transmission efficiency when the system converges, affecting the power increase, thus increasing the number of iterations required. It is worth noting that the significant increase in radiated power is behind the increase in efficiency, which guarantees low power radiation during the low efficiency in the iterative optimization process. This feature effectively reduces the EM wave power radiated to other unknown areas.

Figure~\ref{fig: Fig_spatial_dis} shows the spatial and MT array power density distribution of the EM field radiated by the BS during the iterative process of RF-RBS. The array size is $40 \times 40$, the transmission distance is 2.0~m, and the selected iterations are 1, 2, 5, 10, and 100. The power density is normalized in the figure and does not represent the actual power value. Initially, under the influence of noise, the phase of each element of the BS array is random, resulting in a random and omnidirectional power distribution. The second iteration effectively screens the power density distribution, with most energy radiating towards the MT array. Subsequent iterations focus the energy towards the center of the MT array, though the distribution is not symmetrical. After sufficient iterations, the distribution on the MT array becomes almost symmetric, and the spatial power density distribution shows symmetry between the BS and MT.

The above results show that RF-RBS can optimize the beam through iteration to achieve self-alignment and efficient transmission. Within the effective range, a self-reproducing, resonant beam between BS and MT is established in a finite number of iterations.

\subsection{Power Transmission and Efficiency Optimization}

\begin{figure}
    \centering
    \includegraphics[width=3.5in]{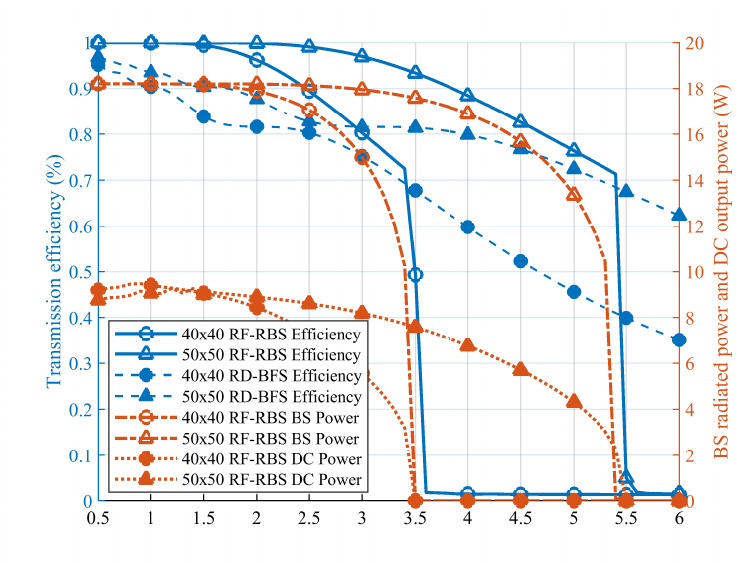}
    \caption{Spatial transmission efficiency $\eta_\text{c}$, BS radiated power $P_\text{BS}$, and DC output power $P_\text{dc}$ versus MT's movement along $z$-axis $\Delta z$.}
    \label{fig: Fig-Efficiency-Power-along-z}
\end{figure}

\begin{figure}
    \centering
    \includegraphics[width=3.5in]{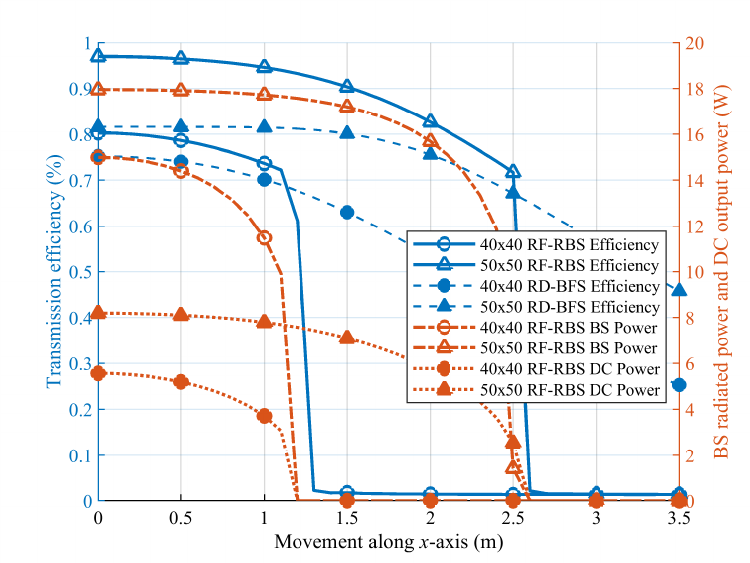}
    \caption{Spatial transmission efficiency $\eta_\text{c}$, BS radiated power $P_\text{BS}$, and DC output power $P_\text{dc}$ versus MT's movement along $x$-axis $\Delta x$.}
    \label{fig: Fig-Efficiency-Power-along-x}
\end{figure}

The RF-RBS achieves self-alignment and efficient transmission by forming a resonant EM beam between BS and MT. This subsection compares and analyzes the power transmission efficiency $\eta_\text{cT}$, BS radiated power $P_\text{BS}$ and DC output power $P_\text{dc}$ of RF-RBS and RD-BFS systems at different array sizes and distances. 

Figure~\ref{fig: Fig-Efficiency-Power-along-z} shows efficiency $\eta_\text{cT}$, BS radiated power $P_\text{BS}$ and DC output power $P_\text{dc}$ of RF-RBS and RD-BFS when MT moves along $z$-aixs with a distance of $\Delta z$. As $\Delta z$ increases, $\eta_\text{cT}$ for both systems initially remains stable but gradually decreases. The RF-RBS consistently outperforms the RD-BFS at medium to short distances, with a maximum efficiency improvement of 16\% when the array size is $40 \times 40$ and the distance is 1.8~m. Additionally, RF-RBS maintains a higher efficiency advantage over the effective distance, while RD-BFS's efficiency declines faster. When the distance exceeds 3.4~m and 5.3~m, respectively, the efficiency of the RF-RBS drops to close to 0, where the transmission efficiency of the RF-RBS can no longer meet the gain-loss condition, and a stable resonance cannot be established. For radiated power $P_\text{BS}$, RF-RBS maintains the highest radiated power at short distances, which is determined by the PAs characteristic. As the distance increases, the radiated power decreases. After reaching the farthest working distance, the radiated power also drops to close to 0, where the radiated power is only the small power radiation under the influence of random noise. 
The DC power $P_\text{dc}$ shows a similar trend to $P_\text{BS}$ and maintains a higher conversion efficiency at shorter distances. The highest conversion efficiency exceeds 50\%. And $P_\text{dc}$ is maintained at W-level within the working distance.
Moreover, systems with more antenna elements achieve higher $\eta_\text{cT}$ values and can support longer transmission distances. By comparing the efficiency changes in Fig.~\ref{fig: Fig-Efficiency-Power-along-z}, with the increase of distance, the attenuation of RF-RBS is smoother, while that of RD-BFS is stepped. 
This is related to the fact that RD-BFS focuses only on the pilot antenna and ignores the other regions of MT, which causes the phase irregularity of the target region \cite{chou2020conformal} and results in non-smooth attenuation of the efficiency.

Similarly, the efficiency $\eta_\text{cT}$, power $P_\text{BS}$ and $P_\text{dc}$ when the MT moves along the $x$-axis with an offset distance of $\Delta x$ is analyzed in Fig.~\ref{fig: Fig-Efficiency-Power-along-x}. The vertical distance along the $z$-axis is set to 3.0~m. When $\Delta x$ is small, both systems exhibit a stable $\eta_\text{cT}$. With the increase of $\Delta x$, the transmission efficiency of both systems decreases, but RF-RBS still outperforms RD-BFS. When $\Delta x$ increases to 1.2~m and 2.5~m, respectively, the efficiency of RF-RBS plummets to 0 because it cannot meet the condition of resonance establishment. The power also shows a similar phenomenon, first maintaining a high level, then gradually decreasing, and after exceeding the effective range, plummeting to close to zero.

The above analysis shows that RF-RBS can form stable EM transmission at different distances and offsets, with a higher transmission efficiency than RD-BFS. In addition, once the transmission efficiency of RF-RBS falls below a specific benchmark, the gain is no longer sufficient to offset the inherent loss. In this case, the resonant beam cannot form between BS and MT, and BS maintains a low radiation power. This characteristic is critical for WPT systems because it minimizes power dispersion to non-target areas, improves overall energy efficiency, and reduces the risk of EM interference. In addition, it can be seen that the minimum transmission efficiency to form resonance is about 70\%, which is independent of the size of the array. From formula \eqref{eq: Gain loss definition}, it can be seen that changing the gain of BS and MT, such as the magnification of PA and the power return ratio of MT, can change the minimum transmission efficiency, thus making a trade-off between transmission efficiency and transmission distance.

\subsection{Power and Phase Distribution}

\begin{figure}
    \centering
    \subfloat[RF-RBS versus RD-BFS with MT at $(0,0,2m)$]{%
    \includegraphics[width=3.5in]{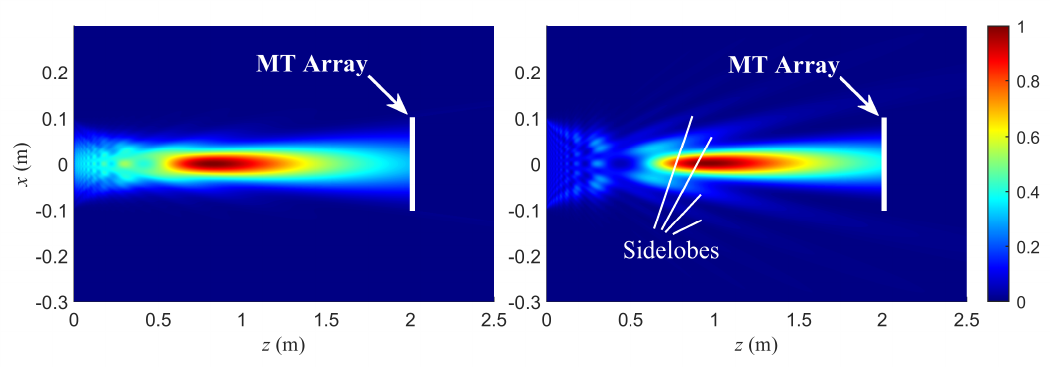}
    } \\
    \vspace{-3mm}
    \subfloat[RF-RBS versus RD-BFS with MT at $(0.5m,0,2m)$]{%
        \includegraphics[width=3.5in]{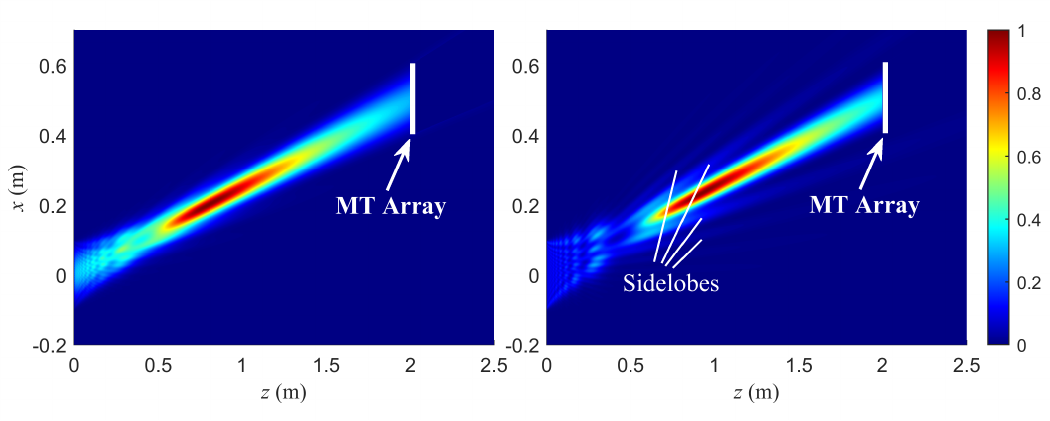}
    }
    \caption{Comparisons of spatial normalized power distribution for RF-RBS and RD-BFS with a $40 \times 40$ array size.}
    \label{fig: Fig-dis-4040}
\end{figure}

\begin{figure}
    \includegraphics[width=3.5in]{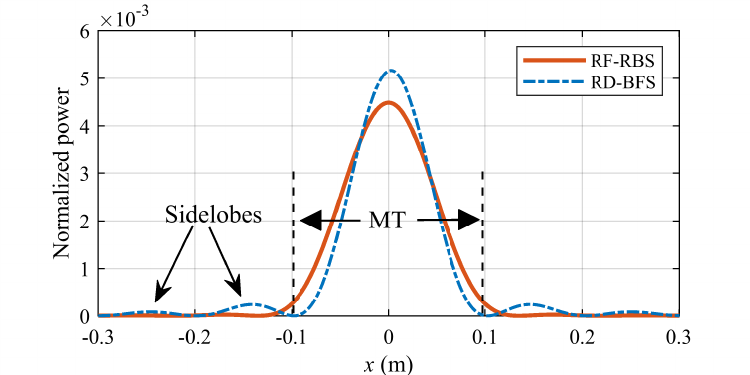}
    \caption{Comparisons of normalized power distribution on MT array for RF-RBS and RD-BFS.}
    \label{fig: Normalized power on MT}
\end{figure}

\begin{figure}
    \centering
    \subfloat[RF-RBS phase distribution with MT at $(0,0,2m)$]{%
    \includegraphics[width=3.3in]{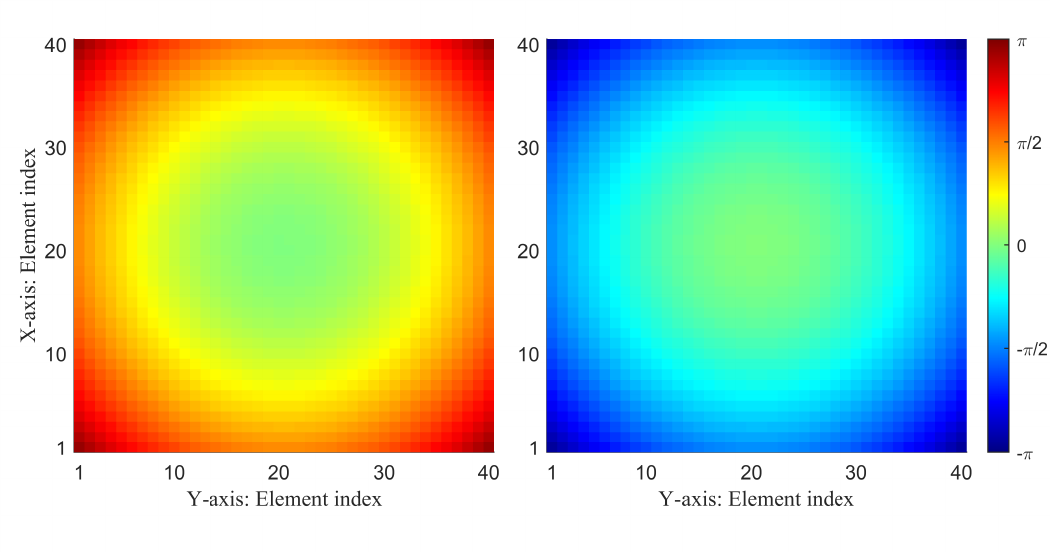}
    } \\
    \vspace{-3mm}
    \subfloat[RF-RBS phase distribution with MT at $(0.5m,0,2m)$]{%
        \includegraphics[width=3.3in]{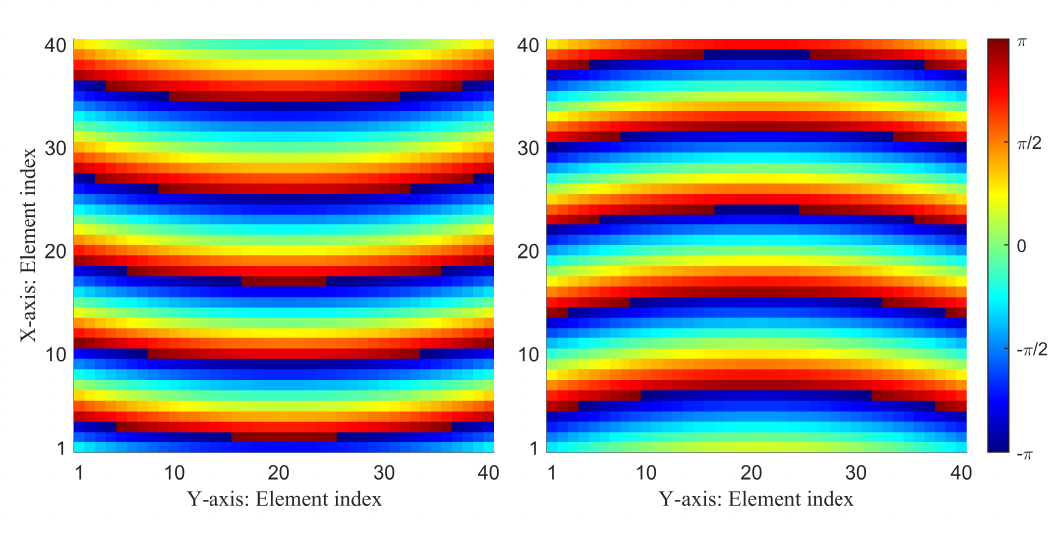}
    }
    \caption{Phase distribution of BS transmission and MT reception.}
    \label{fig: Phase distribution}
\end{figure}

Upon reaching stability, the RF-RBS exhibits a stable power and phase distribution across its arrays, with most energy radiated by the BS directed toward the MT array. This subsection compares and analyzes the spatial power and phase distribution of RF-RBS and RD-BFS.

Figure~\ref{fig: Fig-dis-4040} compares the normalized spatial power density distribution of BS radiation for RF-RBS and RD-BFS with an array size of $40 \times 40$. In Fig.~\ref{fig: Fig-dis-4040}(a), the MT is located at $(0, 0, 2~\text{m})$, and in Fig.~\ref{fig: Fig-dis-4040}(b), the MT is located at $(0.5~\text{m}, 0, 2~\text{m})$. The figures show that RD-BFS has more sidelobes than RF-RBS, preventing part of the energy from radiating to the MT array, resulting in lower efficiency for RD-BFS than RF-RBS. In contrast, RF-RBS has a wider beam, which helps to reduce the generation of sidelobe. Moreover, observing the power distribution, the maximum power point of RF-RBS is closer to the center between the BS and MT, and the power distribution is symmetrically distributed around this point, aligning with resonance characteristics. In contrast, the maximum power point of RD-BFS is closer to the MT. 

Fig.~\ref{fig: Normalized power on MT} shows the normalized power distribution on the MT array when the MT is located at $(0, 0, 2~\text{m})$. The effective receiving range of the MT is from -0.1~m to 0.1~m. As shown in the figure, the power distribution in RF-RBS is smoother, allowing the MT to receive more energy. In contrast, RD-BFS has a higher peak power, but the primary energy is concentrated in the central region of the MT, neglecting the edge regions. RD-BFS also has many sidelobes that fall outside the MT region, resulting in wasted energy, which can lead to additional problems, such as EM interference and safety issues.

Figure~\ref{fig: Phase distribution} illustrates the phase distribution of the BS transmission and the MT reception when resonance is established. When the MT is located at $(0, 0, 2~\text{m})$ in Fig.~\ref{fig: Phase distribution}(a), the phase distribution between the BS and the MT exhibits central symmetry, with isophase lines radiating outward from the center. This symmetrical distribution ensures that the focused signal is directed straight ahead, precisely aligning with the MT's location. When the MT is located at $(0.5~\text{m}, 0, 2~\text{m})$ in Fig.~\ref{fig: Phase distribution}(b), the phase distribution becomes more complex. The isophase lines adjust, radiating outward from the new offset at both the BS and the MT. These lines also demonstrate a tendency to emanate from a specific coordinate point. By fitting the equivalent phase lines and calculating the center coordinates of the circles formed by the isophase lines, the results show that the distance between the antenna's center and the MT's offset along the $x$-axis are nearly equal.

In summary, RF-RBS has fewer sidelobes and a better spatial power distribution than RD-BFS, leading to increased transmission efficiency and improved overall system performance.

\subsection{Communication Performance}

\begin{figure}
    \centering
    \includegraphics[width=3.5in]{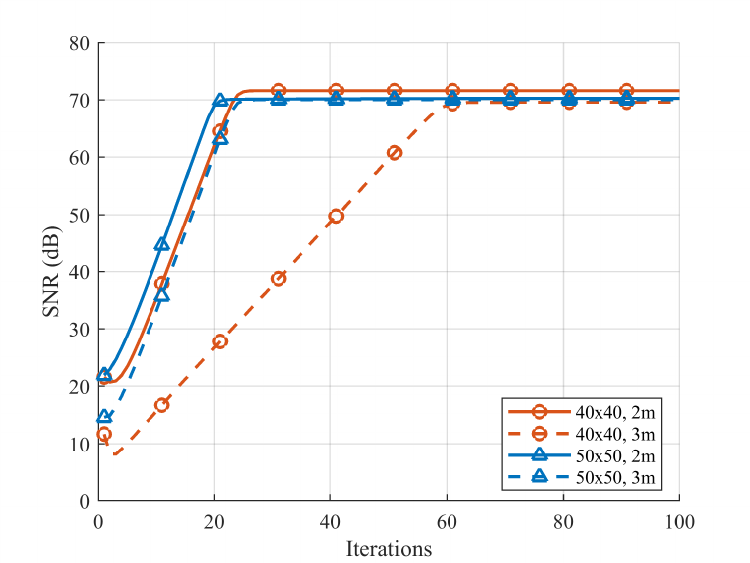}
    \caption{SNR for communication versus the number of iterations.}
    \label{fig: SNR-vs-ite}
\end{figure}

\begin{figure}
    \centering
    \includegraphics[width=3.5in]{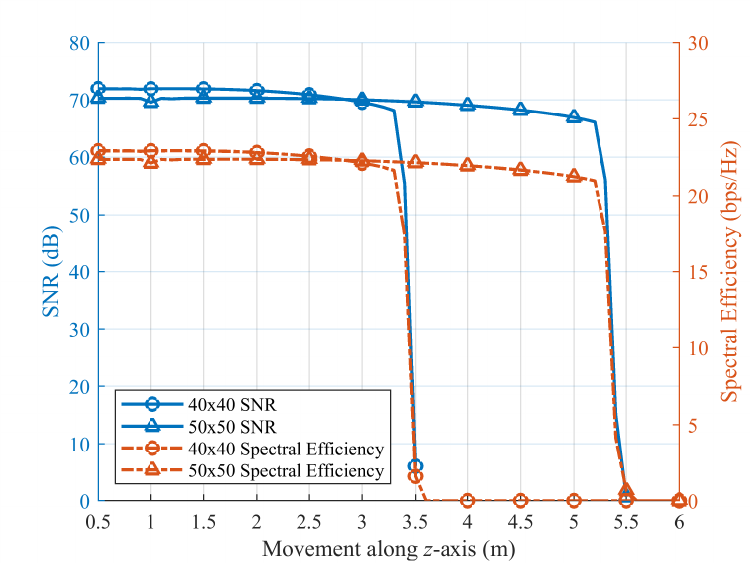}
    \caption{SNR and spectral efficiency for communication versus MT's movement along $z$-axis $\Delta z$.}
    \label{fig: SNR-along-z}
\end{figure}

\begin{figure}
    \centering
    \includegraphics[width=3.5in]{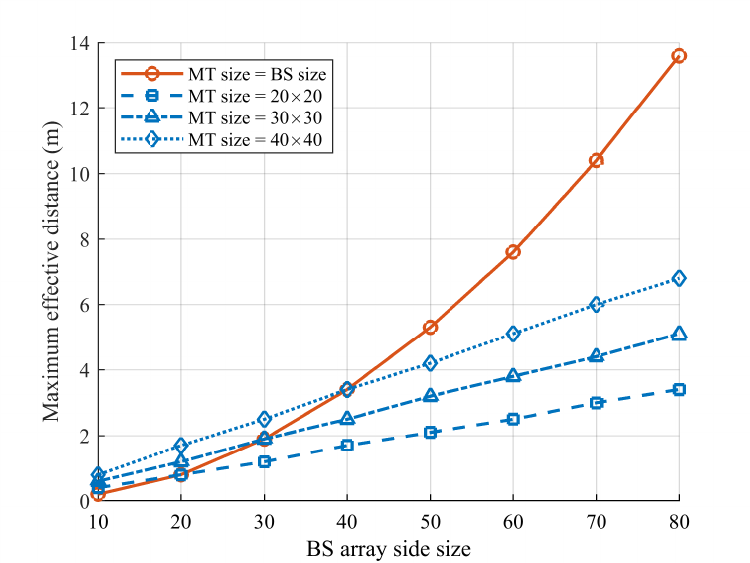}
    \caption{Maximum effective distance versus array sizes of BS and MT.}
    \label{fig: Dis-vs_size}
\end{figure}

Figure~\ref{fig: SNR-vs-ite} shows the changes in SNR for downlink communication with iterations at different array sizes and distances. Before the iteration begins, the signal radiated by the BS is a low-power noise amplified by PAs, and the SNR of the signal received by the MT is infinitesimal. After one iteration, the SNR in the figure improves to more than 10 dB, and the SNR continues to improve with further iterations. After enough iterations, the SNR converges to a stable level.

Figure~\ref{fig: SNR-along-z} shows the SNR and spectral efficiency as the MT moves along the $z$-axis for different array sizes. Due to the resonance effect, the RF-RBS maintains an SNR above 65 dB within effective transmission distances. With a loss factor $\Delta = 3$ dB, the spectral efficiency remains above 21 bps/Hz. For distances less than 2.8~m, the $40 \times 40$ array achieves a higher SNR than the $50 \times 50$ array. This is due to the output power limitations of PAs, which cap the signal power. While the $50 \times 50$ array improves the effective transmission distance, it also introduces more noise, reducing the overall SNR. Despite these limitations, the RF-RBS shows strong potential for significantly enhancing communication channel capacity.

Although current coding techniques and modulation equipment do not fully exploit the calculated maximum capacity, the system's demonstrated feasibility in improving channel performance is promising. In conclusion, the RF-RBS maintains high SNR and spectral efficiency over considerable transmission distances in IoT scenarios. This capability underscores the system's potential to advance communication performance in future wireless power transmission applications.

Additionally, the size of the array affects the maximum effective distance. Figure~\ref{fig: Dis-vs_size} shows the maximum effective distance for different BS and MT array sizes. The horizontal axis of the figure represents the side length of the BS array, whose dimensions are the square of the corresponding values. When the MT is a fixed size, the maximum distance increases linearly with the increase of the BS side length. When the MT and BS increase simultaneously, the maximum distance increases quadratically. The BS and MT have a similar effect on the maximum distance, and their sizes can be adjusted reasonably to meet the needs of different transmission distances. Similarly, a small MT can be used with a large-scale BS for small mobile devices to ensure the transmission distance meets the requirements.

\section{Discussion}
\label{sec:dis}

In the previous sections, the point-to-point system is analyzed. In this section, we discuss multi-target access within the RF-RBS. While an array composed of massive antennas can facilitate multi-directional transmission \cite{xu2014MultiuserMISOBeamforming}, the direction and distance of multiple MTs in RF-RBS can impact transmission loss, making power output to each MT uncontrollable \cite{xiong2019TDMAAdaptiveResonant}. To manage power outputs individually and enhance overall energy utilization, adopting an isolation scheme is essential. In RF-RBS, both Time Division Multiple Access (TDMA) and Frequency Division Multiple Access (FDMA) schemes are feasible.

\subsection{TDMA Scheme}

In RF-RBS, controlled resonance establishment is achievable by toggling the power return on or off in MT. Based on this mechanism, TDMA scheme can be used to achieve multi-target access. By leveraging TDMA, RF-RBS can efficiently manage scenarios involving multiple targets. In this scheme, time is segmented into discrete slots, each designated to a specific MT. During its allocated time slot, the BS verifies the target's activation status and establishes a stable connection, while MTs outside their time slot remain inactive to avoid interference.

MTs within the BS's operational area, denoted as 1, 2, $\dots$, $N_c$, are controlled digitally. Initially, the BS establishes connections with multiple MTs simultaneously, transmitting power and information through omnidirectional radiation, which leads to each MT sharing lower power. Subsequently, the BS allocates time slots to the MTs based on their priority and power needs.

The control in TDMA charging utilizes a digitalized rule comprising frames and slots. Each frame, denoted by $T_f$, is divided into individual time slots allocated to each MT. The allocation follows the rule that the sum of these MTs' active durations must not exceed the frame length:
\begin{equation}
    \sum_{i=0}^{N_c} T_{\text{ON}_i} \leq T_f,
\end{equation}
where $T_{\text{ON}_i}$ represents the active time of the $i$-th MT.

When $N_c$ MTs are present in the work area simultaneously, and they require connection to the BS for power transmission and communication, each MT is allocated time slots within a frame based on its needs and priority. Each MT's time slot includes an additional resonance establishment period ($T_{\text{res}}$) to ensure the formation of a stable resonant connection. Therefore, the average power received by the $i$-th MT is calculated as:
\begin{equation}
    P_{\text{avg},i} = \frac{T_{\text{ON}_i}-T_\text{res}}{T_f} \eta_{c,i}  P_\text{BS} ,
\end{equation}
where $\eta_{c,i}$ is the spatial transmission efficiency from the BS to the $i$-th MT, and $P_\text{BS}$ is the power radiated by the BS.

The primary advantage of TDMA lies in its flexibility in managing power allocation and time distribution. MTs can be prioritized based on their criticality, with high-priority MTs receiving longer time slots and greater power to ensure stable communication and energy transfer. Conversely, lower-priority MTs are allocated shorter slots and reduced power as necessary. Thus, multiple MTs can receive power simultaneously with varying intensities. The BS can generate beams with consistent driving power for each MT, and the driving power can be flexibly adjusted to meet specific requirements.

\subsection{FDMA Scheme}

In RF-RBS, the establishment of resonance between the BS and MTs requires the alignment of their operating frequencies, which are defined by the frequencies of the LO sources in both the BS and MTs. Utilizing this feature, multi-target access can be facilitated by assigning different operating frequencies to various MTs.

In the FDMA scheme, the BS, equipped with massive antennas, is divided into multiple sub-regions, each operating at a distinct frequency band. Typically, in a $2\times2$ sub-array configuration, the BS can be segmented into four sub-regions, each assigned with a different frequency. The MTs within the work area, denoted as 1, 2, $\dots$, $N_c$, tune their LO sources to match the frequency of the sub-region they intend to access. This arrangement allows multiple MTs to establish resonant connections simultaneously, with each operating on a separate frequency.

The total radiated power of the BS is distributed across these frequency bands based on the demand of the MTs:
\begin{equation}
    P_{\text{total}} = \sum_{i=1}^{N_f} P_{\text{BS}_i},
\end{equation}
where $N_f$ represents the number of frequency bands, and $P_{\text{BS}_i}$ is the radiated power allocated to the $i$-th frequency band. By optimizing $P_{\text{BS}_i}$ according to the requirements of the MTs and channel conditions, FDMA can efficiently support multiple MTs.

The primary advantage of FDMA is that it assigns a dedicated frequency channel to each MT, minimizing interference. However, the maximum number of MTs that can be supported is limited by the number of available frequency bands in the BS. The integration of FDMA with the TDMA scheme can further enhance the capability of RF-RBS to support multi-target access.

\section{Conclusion}
\label{sec:con}

This paper demonstrated that RF-RBS, operating in the RF band and within the near-field region, can achieve self-alignment and high-efficiency transmission based on the resonance mechanism. First, we introduced the RF-RBS architecture and the fundamentals of the RAA. Then, we established an analytical model to calculate power flow, resonance establishment, and communication performance. Numerical analysis revealed that the proposed RF-RBS consistently achieves self-alignment without beam control and provides higher transmission efficiency compared to RD-BFS, with improvements of up to 16\%. Additionally, the results indicate the capability to transmit watt-level power and achieve 21 bps/Hz of downlink spectral efficiency in indoor settings. The more concentrated spatial power distribution of the RF-RBS enables high-efficiency transmission over long distances, highlighting the advantages of resonance in WPT and wireless communication applications. Future studies could investigate the system's feasibility and performance under various scenarios and practical implementations.
RF-RBS enables the IoT and mobile devices to operate without cables, presenting application prospects in smart homes, drones, and intelligent manufacturing. However, during implementation, several issues require further investigation, including i) developing multiple safeguards for high-power transmissions, ii) designing high-gain and stable PAs and phase shifters, iii) enhancing communication-perception integration, and iv) integrating Reconfigurable Intelligent Surfaces (RIS) to expand the working range.

%%%%%%%%%%%%%%%%%%%%%%%%%

\ifCLASSOPTIONcaptionsoff
  \newpage
\fi

\bibliographystyle{IEEETran}
\small 

\bibliography{mybib}

\end{document}